\documentclass[12pt,a4paper]{article}
\usepackage[utf8]{inputenc}
\usepackage[english]{babel}
\usepackage[a4paper, total={6in, 10.5in}]{geometry}
\usepackage[OT1]{fontenc}
\usepackage{amsmath}
\usepackage{amsfonts}
\usepackage{amssymb}
\usepackage{amsthm}
\usepackage{tikz}
\usetikzlibrary{positioning, calc}
\usepackage{graphicx}
\usepackage[labelformat=empty]{caption}
\usepackage{float}
\usepackage{bm}
\usepackage{hyperref}
\usepackage[justification=centering]{caption}

\title{On 3-generated 6-transposition groups}
\author{Vsevolod A. Afanasev, Andrey S. Mamontov}
\date{}
\newtheorem*{theorem*}{Theorem}
\newtheorem*{defi*}{Definition}
\newtheorem{lemma}{Lemma}
\newtheorem{obs}{Observation}

\newtheorem{prop}{Proposition}

\newtheorem{notat}{Notation}

\begin{document}
\renewcommand{\refname}{References}
\maketitle

\begin{abstract}
We study $6$-transposition groups, i.e. groups generated by a normal set of involutions $D$, such that the order of the product of any two elements from $D$ does not exceed $6$. We classify most of the groups generated by 3 elements from $D$, two of which commute, and prove they are finite.
\end{abstract}
\section{Introduction}
\let\thefootnote\relax\footnotetext{The work is supported by the Mathematical Center in Akademgorodok under the agreement No. 075-15-2022-281 with the Ministry of Science and Higher Education of the Russian Federation.}
Let $G$ be a group and $n$ be a positive integer. 
If $G$ is generated by a normal set of involutions $D$ such that the order of the product of any two elements from $D$ does not exceed $n$, then $G$ is said to be an \textit{$n$-transposition group} and elements of $D$ are called \textit{$n$-transpositions}. These groups were first considered by Fischer \cite{FischPap} and for the case $n=3$ there exist both a fundamental theory and explicit classification results \cite{3trans,HallCuy}. In particular, $3$-transposition groups are locally finite \cite{HallCuy}, with minimal number of generators for all such finite groups presented in \cite{HallGen3s}.

It is known that the case $n \leq 6$ includes many interesting examples. For instance, the Baby Monster sporadic simple group can be generated by four $2A$-involutions \cite[Theorem 3.1]{SporadGensRels}, that are known to satisfy the $4$-transposition property. The group $B(2,5):2$, an extension of the $2$-generated Burnside group of exponent five (which is yet to be determined to be finite or infinite) by an involution, which inverts generators, can be generated by three $5$-transpositions. Additionally, $2A$-involutions of the Monster sporadic simple group $M$ are $6$-transpositions, with three of them generating~$M$.

The recent development of theories of Majorana algebras and axial algebras has attracted even more attention to the case $n \leq 6$, due to the fact that every Miyamoto group of an axial algebra of Monster type is a $6$-transposition group \cite[Corollary 2.10] {KhaShMC}. We mention the following results on $6$-transposition groups \cite{DecPhD,MSW,KhaShMC} and remark, that they lead to interesting examples and advances in the study of related algebras. 

We aim to classify groups $G=\langle x,y,z \rangle$ generated by three involutions $x,y,z$ from $D$ such that $x$ and $y$ commute, and one of the other two products, say $xz$, is not of order $6$. 

Our main result can be stated as follows:

\begin{theorem*}
If $G=\langle x,y,z \rangle$ is a 6-transposition group, such that $x,\ y,\ z$ lie in its set of $6$-transpositions and $|xy|=2$, $|xz|<6$,
then it is a quotient of one of the following groups: $l^2:D_{12}$ or $l^2:D_{8}$, where $l=4,5,6$; $2^t:D_{2t}$, where $t=5,6$;  $(S_4 \times S_4):2^2$;  $(A_5 \times A_5):2^2$; $PGL(2,9)$; $3^4:(D_8\times S_3)$;  $2 \times (2^s:S_5)$, where $s=4,6$; $k^5:(2^4:D_{10})$, where $k=2,3$; $2^{10}:(2 \times PSL(2,11))$;  $O_2:A_5$ where $|O_2|=2^{10}$; $O_3:D_{20}$ where $|O_3|=3^8$; $2 \times M_{12}$;  $(2.M_{22}):2$; $2\times 2^5:S_6$ or $2 \times 3.S_6$. 
In particular, $G$ is finite.
\end{theorem*}

\textbf{Remark.}  We note that all groups in the statement of the theorem are, in fact, $6$-transposition groups.

Any group from the statement of the theorem can appear as a subgroup of an arbitrary 6-transposition group in the generic case. By generic case we mean that the set of orders of the products of two $D$-elements contains $6$ and either $4$ or $5$, so that one can choose $x,t$ and $z$ with $|xz| \neq 6$ and $|xt|=6$ and using the properties of dihedral groups pick $y \in t^G \subseteq D$ centralizing $x$ from the subgroup $\langle x,t \rangle$. This result forms a stepping stone for further studies regarding $6$-transposition groups.

For the proof the track is twofold: we search for patterns related to $3$-transpositions in order to use the corresponding results for determining the group's structure, and study the involution centralizers to construct explicit presentations and verify them using the GAP system \cite{GAP}. 

The structure of the paper is as follows: we proceed with necessary notation and some elementary observations, then prove the main classification result. After final remarks and discussion, we also provide presentations of the encountered groups in the Appendix.

\section{Notation, preliminaries and results}

Let us fix some notation that is used throughout the rest of the paper:

We write $|X|$ for the order of $X$. When discussing the isomorphism types of groups we use $k$ to denote a cyclic group of order $k$, while $k^n$ is the direct product of $n$ copies of such a group. We use $D_n$ to denote dihedral groups of order $n$, and $M_n$ to denote the Mathieu groups.

We also use $N.H$ for an extension of $N$ by $H$, and $N:H$ if the extension splits. 
We write the elements of $N:H$ as $n.h$, where $n \in N, h \in H$ with multiplication 

$$n_1.h_1\cdot n_2.h_2=(n_1n_2^{h_1}).(h_1h_2)$$ where conjugation means an action of $h_1$ on $n_2$.

Now we highlight the less common notation, which we use extensively:

\begin{notat}
For two elements $g_1, g_2$ of the group $G$ we write $ g_1\sim g_2$ when these elements have equal orders. 
\end{notat}

Recall the definition of a 6-transposition group:

\begin{defi*}\label{6-tr}
A group $G$ is said to be a 6-transposition group if it is generated by a normal set of involutions $D$ such that if $g,h \in D$ then $|gh|\leq 6$.
\end{defi*}

By an \textit{$m$-generated $6$-transposition group} we mean a group, that can be generated by $m$ elements from its set $D$ of $6$-transpositions.

The goal of this paper is to classify all 6-transposition groups with the following Coxeter diagram (assuming $r_1 <6, r_2 \leq 6$):

\begin{figure}[H]
    \centering
    \begin{tikzpicture}[scale=0.6]

    \coordinate[label=left:$x$]  (A) at (0,0);
    \coordinate[label=right:$y$] (B) at (4,0);
    \coordinate[label=above:$z$] (C) at (2,3.464);

    \coordinate[label=below:$2$](c) at ($ (A)!.5!(B) $);
    \coordinate[label=left:$r_1$] (b) at ($ (A)!.5!(C) $);
    \coordinate[label=right:$r_2$](a) at ($ (B)!.5!(C) $);
    
    \draw [line width=1pt] (A) -- (C) -- (B) ;
    \draw[dashed,line width=1pt] (A)--(B);
  \end{tikzpicture}
    \caption{Diagram 1. Coxeter diagram of groups of interest}
    \label{InitDiag}
\end{figure}
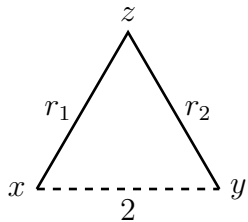
The problem boils down to studying specific homomorphic images of groups satisfying the presentation:
\begin{equation*}\label{presentation}
    C:=\{x,\ y,\ z \mid x^2,\  y^2,\  z^2, (xy)^2, (xz)^{r_1}, (yz)^{r_2} \},\ r_1 <6, r_2 \leq 6. 
\end{equation*}

This restriction, albeit somewhat severe, leaves out
some of the notoriously difficult\ examples of 6-transposition groups such as the aforementioned sporadic simple groups $B$ and $M$, or the group $B(2,5):2$. 

We say that the presentation (or the group) \textit{collapses} or \textit{shrinks} if some of the orders of $xy$, $xz$ or $yz$ are less than $2$, $r_1$ or $r_2$ respectively.

We first state a number of simple yet useful


\textbf{Observations:}

\begin{enumerate}
\item \textit{ Assume $a\in D$ and $i\in G$ is an involution. Then $|ai|$ divides 8, 10 or 12.}  \label{2inv}
    \item \textit{ Assume $a,b \in D$ and $[a,b]=1$. Then $|abc|$ divides 8, 10 or 12 for any $c \in D$.} \label{xyz}
    \item \textit{ Assume $a,b,c \in D$ and $[a,b]=1$. If $|abc|=2k$ is even, then $i=(abc)^k$ is an involution, commuting with $c$ and $ab$. 
    In particular, $\langle a,b,i \rangle \simeq D_{2p}:2$, where $p=|ai|$.
If $|abc|$ is odd, then $ab \in D$.} \label{icommut}

\end{enumerate}

{\it Proof.}
\begin{enumerate}
    \item The square $(ai)^2=aa^i$ is the product of two elements from $D$, thus having order $\leq 6$, so the claim follows.
    \item Corollary of 1.
    \item First assume that  $|abc|$ is even. Then $i$ is the central involution of a dihedral subgroup $\langle ab, c \rangle$. 
Now, $ab$ is in the center of  $\langle a,b,i \rangle$, and in the quotient group two involutions $\overline{a}=\overline{b}$ and $\overline{i}$ generate a dihedral group.
For $|abc|$ odd the claim follows from properties of dihedral groups, since $ab$ is then conjugate to $c \in D$ in $\langle ab,c \rangle$.
$\square $
\end{enumerate}

We additionally provide a useful statement regarding finite $6$-transposition groups.
\begin{lemma}
Let $G$ be a finite $6$-transposition group with a trivial center and $p > 5$ be a prime, then $O_p(G) = 1$.
\end{lemma}

\textit{Proof.}
Suppose that $x\in O_p(G)$ and $a \in D$. Then, as $a$ is an involution, $[x,a]=x^{-1}x^a=a^xa\in O_p(G)$ has order dividing $p$. As $p \ge 7$, $[a,x]=1$ for all $a\in D$ and $x\in O_p(G)$. Hence $O_p(G) \le Z(G)=1$. $\square$


Some of the groups appear for multiple values of $(r_1,r_2)$. In such a case, we record all discovered presentations in the corresponding Appendix entries.
We also use the $p$-core notation, wherever it is not isomorphic to an elementary abelian group. One can use the provided presentations to determine its explicit structure.

Decelle in \cite{DecPhD} studied the $6$-transposition quotients of $C$, additionally requiring that $xy \in D$ (per our notation), and their connections to Majorana theory. Observation \ref{xyz} shows that said condition is rather restrictive.
\section{Proof of the main theorem}

We omit the trivial cases for $r_1,r_2 \in \{1,2\}$ as the only groups that satisfy these conditions are $D_{2n}$ or $2 \times D_{2n}$ for $n \in \{1,6\}$. In all propositions we assume that $G=\langle x,y,z \rangle$, where $x,y,z$ are involutions that satisfy the relations of the diagram \ref{InitDiag}.

\subsection{Case (3,k)} 

By Observation \ref{xyz}, we need the following result for $n \leq 12$:
\begin{prop} Let $n=|xyz|$. If $\langle x,y,z \rangle$ is a quotient of $C$ for $(r_1,r_2)=(3,k),\ k \leq 6,$ then it is a quotient of either
$2 \times A_5$, or $(\frac{n}{2}\times \frac{n}{2}):D_{12}$, where $n$ is even.
\end{prop}

\textit{Proof.} For these cases it follows from \cite[Theorem 6.2.7]{TGT}, that groups, satisfying our diagram are finite for $k\in \{1\ldots 5\}$ (they are the spherical triangle groups) and infinite for $k=6$. Thus, for $k \in \{1\ldots 5\}$ we only need to find the isomorphism types of the groups. For $k=3$ we get that $G \simeq S_4$, where the isomorphism is provided by the map  $x \to (1,2),\ y \to (3,4), z \to (2,3)$.
Similarly, for $k=4$ we get $G \simeq 2 \times S_4$, while for $k=5,\ G \simeq 2 \times A_5$.

We now turn our attention to the case $k=6$.
We will show that the 6-transposition quotients of the group satisfying such a diagram are always finite. 
Set $a=(xzy)^2$ and  $b=(zyx)^2$.
We claim that $[a,b]=1$. 
Indeed,
$$[a,b]=[(xzy)^2,(zyx)^2]=\underline{y}zxyz\underline{x}\cdot \underline{x}yzxy\underline{\underline{z}}\cdot \underline{\underline{xz}}yxzy\cdot zyxz\underline{y}x.$$
By cancelling the underlined letters and changing the double underlined $zxz$ to $xzx$ we obtain
$$zxyzyz\underline{x}y\underline{x}z\underline{x}y\underline{x}zyzyxzx=zx(yz)^5yxzx=zxzxzx=1.$$
Here the last two equalities are due to the fact that $|yz|=6$ and $|xz|=3$.
It follows that $\langle a,b \rangle$ is the direct product of two cyclic groups of order $|a|$.
Since $a^x=b$, $a^z=ab^{-1}$, and $a^y=b^{-1}$, we conclude that $\langle a,b \rangle$ is normal in $\langle x,y,z \rangle$.

We now consider the corresponding quotient group, using same notation for the images of elements. Note that $1=a^z=xz^{yz}$, so $x \in \langle y,z \rangle$.
Hence the quotient group is isomorphic to $D_{12}$. The proposition follows for even $n$.

For odd $n$, the elements $a$ and $b$ generate the group $H=\langle yzx,xyz\rangle$ which, by previous reasoning, is abelian and normal. In the quotient $G/H$ we get that $xy=z$ and in particular $[x,z]=1$, but since the order of $xz$ should divide 3, implying that $x=z$ (in $G/H$) this means that $y=1$ and thus $y \in H$. Then $xz \in H$ and since $H$ is abelian we get that $[y,xz]=1=yzxyxz=(yz)^2$, which means that the group shrinks. 

Finally note that the Coxeter group $2 \times S_4$ is isomorphic to $(2\times 2):D_{12}$, and so it is covered by the statement of this proposition. $\square$

\subsection{Case (4,4)}
In this case we will restrict the order of $xy^z$ to obtain the following

\begin{prop}
Assume $(r_1,r_2)=(4,4)$ and set $k$ to be the order of $xy^z$. Then $G$ is an extension of a direct product $k\times k$ of two cyclic groups of order $k$ by $D_8$. 
\end{prop}

\textit{Proof.} Consider the group $K=\langle x,x^z,y,y^z \rangle$. The index of $K$ in $G$ is $\leq 2$. Let us compute the relations between pairs of involutions generating $K$. Conjugating the relation $[x,y]=1$ we obtain $[x^z,y^z]=1$. The order of $xx^z=(xz)^2$ is $2$, hence $[x,x^z]=1$. Similarly $[y,y^z]=1$. Overall, we get that $[\langle x^z,y \rangle, \langle x,y^z \rangle]=1$. 

Note that $|x^zy|=|xy^z|=k$. Therefore $\langle x^z,y \rangle \simeq D_{2k}$, in other words it is an extension of a cyclic group $\langle x^zy \rangle$ of order $k$ by an involution.
It is clear now that $L=\langle x^zy, xy^z \rangle$ is a direct product $k \times k$ of two cyclic groups of order $k$. Obviously, $L$ is normal in $K$. Conjugating the generators of $L$ by $z$ permutes them, so $L$ is also normal in $G$. In the quotient $G/L$ the image $\overline{y}=\overline{x^z}\in \overline{\langle x,z \rangle}$. Therefore $G/L \simeq D_8$. The proposition is proved. $\square $

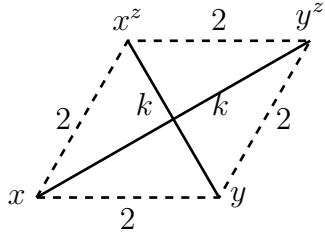
\begin{figure}[ht]
    \centering
    \begin{tikzpicture}[scale=0.6]

    \coordinate[label=left:$x$]  (A) at (0,0);
    \coordinate[label=right:$y$] (B) at (4,0);
    \coordinate[label=above:$x^z$] (C) at (2,3.464);
    \coordinate[label=above:$y^z$] (D) at (6,3.464);

    \coordinate[label=below:$2$](c) at ($ (A)!.5!(B) $);
    \coordinate[label=left:$2$] (b) at ($ (A)!.5!(C) $);
    \coordinate[label=left:$k$](a) at ($ (B)!.6!(C) $);
    \coordinate[label=above:$2$](a) at ($ (D)!.5!(C) $);
    \coordinate[label=right:$k$](a) at ($ (A)!.6!(D) $);
    \coordinate[label=right:$2$](a) at ($ (D)!.5!(B) $);
    \draw [dashed,line width=1pt] (A) -- (C) -- (D)--(B)-- cycle;
    \draw[line width=1pt] (A)--(D);
     \draw[line width=1pt] (B)--(C);
  \end{tikzpicture}
    \caption{Diagram 2. Diagram for subgroup $K$.}
    \label{CoxKGroup}
\end{figure}

\subsection{Case (4,*)}

The goal of the section is to complete the case $r_1=4$ and to prove the following:

\begin{prop}
Assume $r_1=4$, $r_2 \in \{5,6\}$. Then $G$ is the quotient of one of the following groups:
$2^{r_2}:D_{2r_2}$, $2\times (2^s:S_5)$, where $s=4,6$, $3^4:(D_8 \times S_3)$, $(S_4 \times S_4):2^2$, $(A_5\times A_5):2^2$, or $PGL(2,9)$.
\end{prop}

We will begin by showing that

\begin{lemma}\label{5-gend}
$\langle x^G \rangle =\langle x, x^z, x^{zy}, x^{zyz}, x^{(zy)^2}, x^{(zy)^2z} \rangle$.
\end{lemma}
For $r_2=5$, the first $5$ generators suffice, since $(zy)^2z=yzyzy$ and $x^{(zy)^2z}=x^{(zy)^2}$.

\textit{Proof.}
Set $L= \langle y,z\rangle \cong D_{2r_2}$. Then as $y$ and $x$ commute, $x^{L}$ has order dividing $r_2$ and is the set $\{x,x^2,x^{zy},x^{zyz}, x^{(zy)^2}, x^{(zy)^2z}\}$. Hence $N=\langle x,x^2,x^{zy},x^{zyz}, x^{(zy)^2}, x^{(zy)^2z}\rangle$ is normalized by $L$ and, as $x \in N$, $N$ is normalized by $G$. This proves the claim.
$\square$

We can additionally show that 

\begin{lemma}\label{prodorderz}
The product orders for generators of $\langle x^G \rangle$ can be expressed via orders of two elements 
 of $D$: $x\cdot x^z$, $x\cdot x^{zyz}$ and $x \cdot x^{(zy)^2z}$. 
\end{lemma}
\textit{Proof.}
We see that 
$xx^{z}=(xz)^2$,
$xx^{zyz}=x(zyz)x(zyz)=(xy^z)^2$,
$xx^{(zy)^2z}=(xz^{yz})^2$.

For other products the statement can be checked via conjugation of the previous three formulas. For instance,
$x^zx^{zyz}=(xx^{zy})^z=((xz^x)^y)^z$. $\square$

\textbf{Remark.} For $r_2=5$, only the first two formulas suffice since $x^{(yz)^2z}=x^{(yz)^2}$.

We have that for even values of $|xy^z|$ and $|xz^{yz}|$ the relations between generators of $\langle x^G \rangle $, resemble those of 3-transposition groups, a property which will allow us to determine the structure of our 6-transposition groups.

We will now proceed with the case-by-case analysis for values of $|xy^z|$ and $|xz^{yz}|$.

\begin{lemma}
For $|xy^z|=4$
the group $G$ is a quotient of $2^{r_2}:D_{2r_2}$.
\end{lemma}
\textit{Proof.}
For case $|xz^{yz}|=4$ it can be seen that $\langle x^G \rangle$ is elementary abelian of order $2^{r_2}$. Note that the index of $G:\langle x^G \rangle$ does not exceed $2r_2$ and we have an action of $\langle y,z \rangle$ on $\langle x^G \rangle$.
We provide an isomorphism, describing this action. 
Treat $\langle x^G \rangle$ as a vector space over  $\mathbb{F}_2$ with basis $e_i$ and the generating involutions of $D_{2r_2}$ as permutations, that correspond to reflections of a polygon, while fixing a point. 
For $r_2=5$ take 
$$x=e_1.; y=.(2,5)(3,4); z=.(1,3)(4,5);$$
so that conjugation is treated as permutation of indices of basis vectors.
It is straightforward now to check the relations, and confirm the isomorphism.

For $r_2=6$ the reasoning is entirely similar. 

Moving on to the case $|xz^{yz}|=6$. For $|yz|=5$ the diagram is similar to the previous one. So we further assume $|yz|=6$.
By Lemma \ref{prodorderz}, the diagram looks like three lines in the Fisher space, as shown below.

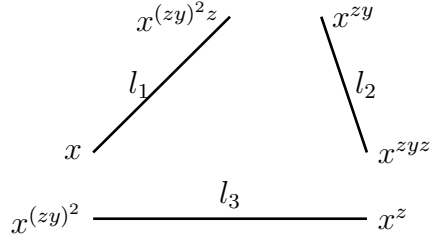
\begin{figure}[H]
    \centering
     \begin{tikzpicture}[scale=0.6]

    \coordinate[label=left:$x$]  (A) at (0,0);
    \coordinate[label=right:$x^{zyz}$] (B) at (6,0);
    \coordinate[label=right:$x^z$] (C) at (6,-1.464);
    \coordinate[label=right:$x^{zy}$] (D) at (5,3);
    \coordinate[label=left:$x^{(zy)^2}$] (F) at (0,-1.464);
    \coordinate[label=left:$x^{(zy)^2z}$] (G) at (3,3);

    \coordinate[label=above:$l_3$](c) at ($ (C)!.5!(F) $);
    \coordinate[label=right:$l_2$] (b) at ($ (D)!.5!(B) $);
    \coordinate[label=left:$l_1$](a) at ($ (A)!.5!(G) $);
     \draw[line width=1pt] (F)--(C);
     \draw[line width=1pt] (D)--(B);
    
     \draw[line width=1pt] (A)--(G);
    
  \end{tikzpicture}
    \caption{Diagram 3. Diagram for the case $|xy^z|=4,\ |xz^{yz}|=6$.
    }
\end{figure}

Therefore $G$ is a 6-transposition quotient of $K \simeq S_3^3:D_{12}$. Observe that the generating elements $x$, $y$, and $z$ in $K$ act on the points of the Fisher diagram for $S_3^3$ in the following way:\newline
$x$ fixes lines $l_2$ and $l_3$, and permutes two of the three points on $l_1$;\newline
$y$ interchanges $l_2$ and $l_3$ and centralizes points of $l_1$;\newline
$z$ interchanges $l_1$ and $l_3$, and permutes two points on $l_2$.\newline
It follows that $xyz$ transitively permutes $9$ points of the Fisher space. By Observation~\ref{xyz}, $G$ should be a quotient of $K$ by $\langle (xyz)^6 \rangle$. Points on each line of the Fisher diagram will merge correspondingly into a single point, and so $G$ is a quotient of $L\simeq 2^3:D_{12}$. In $L$ we have $(xyz)^3,(yz)^3 \in Z(L)$ and $|xz^{yz}|=2$, so $L$ can be also recognized as $2^2 \times S_4$, and it is a quotient of the group from the previous subcase.

Finally we will see that for $|xz^{yz}|=5$ the group shrinks. When $|yz|=5$ this is obvious, since 
$xz^{yz}=xy^{zy}\sim x^yy^z=xy^z,$
however we assume that the last element has order 4.
So we consider $|yz|=6$.
Then $\langle x^G \rangle$ should be the quotient of $D_{10}^3$. Since $|xz^{yz}|=5$ we have $z \in x^G$, so it can permute the lines corresponding to $D_{10}$ only if the group shrinks. $\square$

\begin{lemma}\label{starSix}
For $|xy^z|=6$ the group $G$ is a quotient of $2 \times (2^4:S_5)$, $2\times (2^6:S_5)$, $3^4:(D_8 \times S_3)$ or $(S_4 \times S_4):(2^2)$.
\end{lemma}
The relations between the elements are shown in the following diagram.

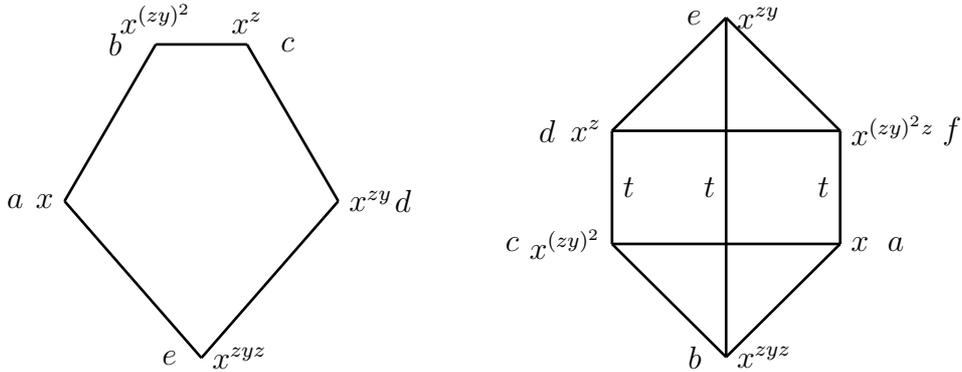
\begin{figure}[H]\label{CoxDiag3}
    \centering
    \begin{tikzpicture}[scale=0.6]

    \coordinate[label=left:$x$]  (A) at (0,0);
    \coordinate[label=right:$x^{zyz}$] (B) at (3,-3.464);
    \coordinate[label=above:$x^z$] (C) at (4,3.464);
    \coordinate[label=right:$x^{zy}$] (D) at (6,0);
    \coordinate[label=above:$x^{(zy)^2}$] (F) at (2,3.464);
    \coordinate[label=right:$a$] (AA) at (-1.5,0);
     \coordinate[label=left:$b$] (FF) at (1.5,3.464);
     \coordinate[label=right:$c$] (CC) at (4.5,3.464);
    \coordinate[label=right:$d$] (DD) at (7,0);
    \coordinate[label=left:$e$] (BB) at (2.7,-3.464);

    \draw [line width=1pt] (A) -- (B);
    \draw[line width=1pt] (A)--(F);
     \draw[line width=1pt] (D)--(C);
     \draw[line width=1pt] (F)--(C);
     \draw[line width=1pt] (D)--(B);
  \end{tikzpicture}
  \qquad
    \begin{tikzpicture}[scale=0.6]

    \coordinate[label=right:$x$]  (A) at (5,0);
    \coordinate[label=right:$x^{zyz}$] (B) at (2.5,-2.5);
    \coordinate[label=left:$x^z$] (C) at (0,2.5);
    \coordinate[label=right:$x^{zy}$] (D) at (2.5,5);
    \coordinate[label=left:$x^{(zy)^2}$] (F) at (0,0);
    \coordinate[label=right:$x^{(zy)^2z}$] (G) at (5,2.5);

     \coordinate[label=right:$a$] (AA) at (5.8,0);
     \coordinate[label=left:$c$] (FF) at (-1.8,0);
     \coordinate[label=left:$d$] (CC) at (-1,2.5);
    \coordinate[label=right:$f$] (GG) at (7,2.5);
    \coordinate[label=left:$e$] (DD) at (2.2,5); 
    \coordinate[label=left:$b$] (BB) at (2.2,-2.5);

    \coordinate[label=left:$t$](c) at ($ (A)!.5!(G) $);
    \coordinate[label=right:$t$](cc) at ($ (C)!.5!(F) $);
     \coordinate[label=left:$t$](ccc) at ($ (B)!.5!(D) $);
    \draw [line width=1pt] (A) -- (B);
    \draw[line width=1pt] (A)--(F);
     \draw[line width=1pt] (D)--(C);

      \draw[line width=1pt] (C)--(G);
       \draw[line width=1pt] (D)--(G);

    \draw [line width=1pt] (A) -- (G);
   \draw [line width=1pt] (C) -- (F);
    \draw [line width=1pt] (B) -- (D);
     
    \draw[line width=1pt] (B)--(F);
     
  \end{tikzpicture}
    \caption{Diagram 4. Diagram for subgroups $\langle x^G \rangle$ when $|xy^z|=6$,\newline for $r_2=5$ (left) $r_2=6$ (right). $t=|xz^{yz}|$.
    }
    \label{AffCoxDiagms}

\end{figure}
\textit{Proof.}
First assume $|yz|=5$. 

By \cite[Appendix, 14a]{HallSoicher}, if the group $\langle a,b,c,d,e\rangle$ has generators satisfying the diagram of type $\widetilde{A}_4$ and if element $l=a^{bcd}e$ has order $\leq 3$ then the resulting group is a $3$-transposition group.
More precisely, $\langle a,b,c,d,e\rangle \simeq~ k^4:S_5$, where $k=|l|$.
To test this condition, we reduce the aforementioned element $l$ by setting $a:=x$ in our previous notation and going clockwise over Diagram \ref{AffCoxDiagms}, (that is, $b:=x^{(zy)^2}$, etc.)
$$a^{bcd}e=x^{x^{(zy)^2}x^zx^{zy}}x^{zyz}.$$

We first note that $a^b=x^{(zy)^2}xx^{(zy)^2}=$ $y(zyzx)^3zyzy=y(xzyz)^3zyzy=y(xy^z)^2yx$. Conjugating this by $c$ we get 
$zxzyxy^zxy^zyxzxz$. Since $|xz|=4$, we push both $x$ to the left and right and get
$$xzxzyy^zxy^zyzxzx.$$ Using $zyzyz=yzyzy$ we obtain $$xzxyy^zyxyy^zyxzx=xzy(xy^z)^2yxzx.$$
We compute $a^{bc}d=a^{bc}x^{zy}=x(zyxy^zxzy)^2x$ (pushing the last $x$ to the right).

Then $a^{bc}de=(zyx)^4y^zyx$ and adding $d$ to the left we simply get $yzxzy(zyx)^4y^zyx$. We then note that we can rewrite this element as follows:
$$yzxzy(zyx)^4y^zyx\sim zyzxzxzy(zyx)^4=(zyx)zxzzy(zyx)^4=(zyx)^6.$$

By Observation \ref{2inv}, $k=|(xyz)^6| \leq 5$.

Moving on to the group $G$ we first note, looking at the diagram, that\newline
$(xyz)^5=x\cdot x^{zy}\cdot x^{(zy)^2}\cdot x^{(zy)^3}\cdot x^{(zy)^4}=adbec \in \langle x^G \rangle,$\newline
and for $a\simeq (1,2)$, $b\simeq (2,3)$, $c\simeq (3,4)$, $d\simeq (4,5)$, $e\simeq (1,5)$, this equals to $(1,4)(2,5,3)$ and has order $6$.
We conclude by Observation \ref{2inv}, that $zy$ must lie in $\langle x^G \rangle$, and $(xyz)^{12}=1$, so the order of $xyz$ cannot be $8,5$, or $10$, since $S_5$ has no corresponding quotients. Therefore, $|G:\langle x^G \rangle|=2$, and $G\simeq 2\times (2^4:S_5)$.
Using the representation $x\simeq 1 \times e_2.(1,2)$, $y\simeq 1 \times e_1.(1,2)(4,5)$, $z \simeq 1 \times e_2.(1,5)(2,3)$ it is straightforward now to check the relations and the 6-transposition property.

\vspace{5pt}

{\it Case $|yz|=6$ and $|xz^{yz}|=4$.}\newline
Looking at Diagram \ref{AffCoxDiagms} (with labeled edges removed, since $t=2$), we observe that \newline
1) $z$ permutes the triangles;\newline
2) $y$ fixes one vertex in each triangle, namely $x$ and $x^{(zy)^2z}$, while permuting two other.\newline
The order of the product of vertex elements 
$x \cdot x^{zyz}\cdot x^{(zy)^2} = (xy^zy)^3$
is not divisible by~$3$ by Observation \ref{xyz}. 
Denoting the order of this product by $2k$, where $k\in \{1,2,5\}$, it follows that elements in each triangle generate a quotient of $(k\times k):(3\times 3):2$, with the 3-transposition group on top.
Let us again denote vertices of triangles as per the second part of Diagram \ref{AffCoxDiagms}.
Then $[\langle a,b,c \rangle,\langle d,e,f \rangle]=1$. We also have
\begin{equation} \label{ll1}
(xyz)^6 = x\cdot x^{zy}\cdot x^{(zy)^2}\cdot x^{(zy)^3}\cdot x^{(zy)^4}\cdot x^{(zy)^5}=aecfbd=acb \cdot efd,
\end{equation}
which is a product of two commuting elements of order $2k$. We deduce that $k\neq 5$.

First assume $k=1$. Then $(xyz)^{12}=1$, and $\langle x^G \rangle$ is a 3-transposition quotient of $[(3\times 3):2] \times [(3\times 3):2]$.
With $D_{12}$ on top of that, we obtain that $G$ is a quotient of $3^4:(D_8 \times S_3)$.

Finally assume $k=2$. Then $(xyz)^8=1$ and $\langle x^G \rangle$ is a quotient of the direct product of two groups, isomorphic to $(3 \times A_4):2$, where the top involution acts as it would on the corresponding subgroups in $S_3$ and $S_4$  From (\ref{ll1}) we also have $(xyz)^2 \in \langle x^G \rangle$, hence $(yz)^2 \in \langle x^G \rangle$, and so $G$ can be recorded as a quotient of $(S_4\times S_4):2^2$.

\vspace{5pt}

{\it Case $|yz|=6$ and $|xz^{yz}|=5$.}\newline
In this case $G$ shrinks. This can be seen by looking at the following diagram:
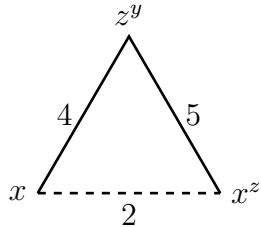
\begin{figure}[H]
    \centering
    \begin{tikzpicture}[scale=0.6]

    \coordinate[label=left:$x$]  (A) at (0,0);
    \coordinate[label=right:$x^z$] (B) at (4,0);
    \coordinate[label=above:$z^y$] (C) at (2,3.464);

    \coordinate[label=below:$2$](c) at ($ (A)!.5!(B) $);
    \coordinate[label=left:$4$] (b) at ($ (A)!.5!(C) $);
    \coordinate[label=right:$5$](a) at ($ (B)!.5!(C) $);
   
   \draw [line width=1pt] (A) -- (C) -- (B) ;
    \draw[dashed,line width=1pt] (A)--(B);
  \end{tikzpicture}
     \caption{Diagram 5. Diagram for the subgroup, case $|yz|=6,\ |xz^{yz}|=5$.}
    \label{CoxDiag}
\end{figure}

\noindent In this subgroup
$x\cdot x^{z(z^y)}\sim xx^{zyz}$ is of order $3$, but while analysing the quotients of the group we previously obtained, it can be seen that this subgroup shrinks, thus shrinking the entire  $G$.

{\it Case $|yz|=6$ and $|xz^{yz}|=6$.}\newline
From Diagram \ref{AffCoxDiagms}, it is clear that $\langle x^G \rangle$ has a quotient, isomorphic to $S_5$, while $\langle a,c,d,e \rangle \simeq S_5$.
Note that $y,z$ act on the generators of this group as the following permutations:\newline
$y=(c,b)(d,e),$ 
$z=(a,d)(c,f)(e,b).$ \newline
We consider the subgroup $\langle x, z^y, y^z \rangle$. Note that the product of generators $xy^zz^y=x(zy)^3$ should be of order $6$, 
since $(x(zy)^3)^2=xx^{(zy)^3}\sim xx^{(zy)^2z}$ has order 3.

The last product of Lemma 2 for this group looks as follows:  \newline
$z^y\cdot x^{z^yx}\sim zx^{zyxy}\sim zx^{zx}=xz$ and has order 4.
These two observations and previous cases imply that $|z^yy^{zx}|=4$.
By Lemma 3 we then have $\langle x,z^y,y^z \rangle \simeq 2\times (2^4:D_{12})$.

It is clear that $a,b,e,f \in \langle x,z^y,y^z \rangle$.
Direct computations in this group show that the following 3-transposition relations hold:
$(a^{fe}b)^2=(afeb)^3=1$. \newline
Conjugating by $y=(c,b)(d,e)$, we also obtain relations
$(a^{fd}c)^2=(afdc)^3=1$. \newline
These relations for the diagram imply that $\langle x^G \rangle \simeq 2^5:S_5$.
Going back to $G$ as before we have $(yz)^2 \in \langle x^G \rangle$. 
Therefore $G$ is $2\times 2^6:S_5$. $\square$
\begin{lemma}
For $|xy^z|=5$ the group $G$ is a quotient of $(A_5\times A_5):2^2$, or $PGL(2,9)$.
\end{lemma}

\textit{Proof.}
For $|yz|=5$ we notice that the subgroup $\langle xy,\ xz \rangle$ satisfies the well known presentation of $A_6 \simeq \langle a,b \mid a^2,b^4,(ab)^5,(ab^2)^5\rangle$ (in this presentation $a\simeq (1,2)(3,4)$ and $b\simeq (1,2,3,5)(4,6)$) and is clearly normal in $G$. The quotient order is equal to $2$ and thus we only need to determine the extension type of this involution to determine the structure of $G$. We show that this group is isomorphic to $PGL(2,9)$ by explicitly presenting the matrices that satisfy the relations for $x,y,z$:
$$x\simeq\left( \begin{array}{cc}
    1+2\alpha & 2+\alpha \\
    2+2\alpha & 2+\alpha 
\end{array}\right),\; y \simeq\left( \begin{array}{cc}
    1& 2\\
    0 & 2 
\end{array}\right)\; z\simeq\left( \begin{array}{cc}
    2\alpha & 2\alpha \\
    2+2\alpha & 2+\alpha 
\end{array}\right),
$$
where $\alpha$ is a root of an irreducible polynomial of degree $3$ over $\mathbb{F}_3$. 

Now set $|yz|=6$. 
Then $G$ is a homomorphic image of 
$$G(i,j)=\langle x,y,z \mid K \cup \{(xy^z)^5, (xyz)^{i}, (xz^{yz})^j\} \rangle,$$
where $K$ is the set of Coxeter relations, $i \in \{8,10,12\}, j \in \{4,5,6\}$. 
We claim that the only possibility is $G\simeq G(12,6)\simeq (A_5 \times A_5): 2^2$.
Using coset enumeration we confirm the claim in terms of group orders. 
Meanwhile the following permutations:\newline
$x=(1,2)(3,4)(5,6)(7,8)(9,10)(11,12),\\
y=(1,2)(3,4)(5,12)(6,11)(7,10)(8,9),\\
z=(1,12)(2,7)(3,6)(4,5)(8,10)(9,11)\\$
satisfy the relations and generate a subgroup of the desired isomorphism type. $\square$

{\bf Proof of Proposition 3} is provided by Lemmas 2-6. $\square$

\subsection{Case (5,5)}

Moving along with our case-by-case analysis, this section's goal is to prove

\begin{prop}\label{5n5}
Assume $(r_1,r_2)=(5,5)$. Then $G$ is the quotient of one of the following groups: $PGL(2,9)$, $k^5:(2^4:D_{10})$, where $k=2,3$; $2^{10}:(2 \times PSL(2,11))$ or $O_2(G):A_5$, where $|O_2(G)|=2^{10}$.
\end{prop}

\setcounter{obs}{3}
\begin{obs}\label{WholeGroupObs}
If $r_1=5$ then $z \in \langle x,x^z\rangle $. If $r_2=5$ then $ z \in \langle y, y^z \rangle$.
\end{obs}

We introduce the notation that helps structure further proofs:
\begin{notat}
     We write $G(i_1, \ldots i_n):= \langle x, y, z \mid K \cup \{ w_1^{i_1},\ldots, w_n^{i_n} \} \rangle $, where $w_i=w_i(x, y, z)$ are fixed, and $K$ are the Coxeter relations used throughout the section.
\end{notat}

\textit{Proof of Proposition \ref{5n5}.}

Set $w_1=xy^z,w_2=xyz\text{ and } w_3=xz^{(xyz)^2}$.
We iterate in GAP \cite{GAP}, to calculate the order of $G(i_1, i_2, i_3)$ for $i_1,i_3 \in \{4,5,6\}$ and $i_2 \in \{8,10,12 \}$.
This leaves us with just a few non-shrinking presentations. We identify their 6-transposition quotients, recording additional relations for nontrivial cases.

Groups with $i_1=4$ were already considered since in this case $G=\langle x,y^z,y\rangle$.

The group $G=G(5,12)$ is isomorphic to $\langle i^G \rangle: (2^4:D_{10})$. In this case $|\langle i^G \rangle|=3^5$, where $i=(xyz)^4$.
By the Schur-Zassenhaus theorem the extension splits, and so $G \simeq \langle i^G \rangle: (2^4:D_{10})$. 
Moreover, the elements $i,i^x,i^{xz},i^{xzx}$ and $i^{xy^z}$ generate $\langle i^G \rangle$ and commute pairwise. 

We find that $G(5,8) \simeq \langle i^G \rangle: (2^4:D_{10})$ for the same word $i$.

In turn, $G(6,10,4)=\langle j^G \rangle : G(6,10,2)$ with $\langle j^G \rangle$ having order $2^{10}$ and being elementary abelian for $j=(xz^{(xyz)^2})^2$, while $G(6,10,2)\simeq 2 \times PSL(2,11)$.
For $(G,10,6)/\langle(xyy^z)^6\rangle$  we again obtain $2 \times PSL(2,11)$.

Finally the group $G(6,12,4)$ is finite of order $2^{12}\cdot 3 \cdot 5$ and it satisfies the $6$-~transposition property. Further research shows, that it has a normal subgroup of order $2^{10}$, which is a normal closure of $\langle (xyz)^3 \rangle$ and it is equal to $O_2(G)$, with the quotient group isomorphic to $A_5$.
$\square$

\subsection{Case (5,6)}
\begin{prop}
Assume $(r_1,r_2)=(5,6)$. Then $G$ is a quotient of one of the following groups: $2 \times (2^4:S_5)$, $2 \times 2^5:S_6$,
$O_3(G):D_{20}$, where $|O_3(G)|=3^8$, $2 \times M_{12},\ (2.M_{22}):~2,\  2 \times 3.S_6$ or one of the groups from Proposition \ref{5n5} (sans $PGL(2,9)$).

\end{prop}

\textit{Proof.}
Using Observation \ref{WholeGroupObs}, we can instantly deal with two cases:
letting $|xy^z|=|x^zy|=4$ the group $\langle x, y, x^z \rangle=G$ satisfies the $(2,4,5)$ diagram. The corresponding relators look as follows: $(yx^zxx^z)\sim yz$ is of order $6$, while $(xyx^z)^6$, per reasoning in Lemma \ref{starSix}, can have order $1$ or $2$. If we have it set to $2$ we get $G \simeq 2 \times (2^4:S_5)$.

Next, if we suppose $|xy^z|=|x^zy|=5$ the group  $\langle x, y, x^z \rangle$ satisfies the $(2,5,5)$ diagram and is again equal to $G$. It is either $2^{10}:(2 \times PSL(2,11))$ or $O_2(G):A_5$.

Therefore we set $|xy^z|=6$ for the rest of this section. 

We will start with a construction, that allows us to obtain a centralizer element for $x$. Consider the subgroup $\langle x, z,w \rangle$, where $w=z^y$. We have that $|xz|=5$, $|zw|=|zyzy|=3$ and $|xw|=|xyzy|=|xz|=5$. Thus all the generator products have odd order. This implies that if we take any pair of these generators and consider the dihedral group generated by that pair, then these elements will be conjugate in said group, in fact, we can write the conjugating elements as follows:
$$x^{(wx)^2}=w=z^{wz}.$$ 
Since $x$ and $z$ are also conjugate we get that 
$$x^{(wx)^2}=x^{(zx)^2(wz)},$$ 
in other words, the element $(zx)^2z^yz(xz^y)^2 \in C_G(x)$.
Rewriting this, we get $$(zx)^2(yz)^2yxzyxyzy=(zx)^2(yz)^2yxzxzy.$$ 
Since $y \in C_G(x)$ we can remove it from the right to get 
$(zx)^2(yz)^2y(xz)^2 \in C_G(x)$, in other words $y^{zy(xz)^2} \in C_G(x)$.

We will denote the obtained centralizer element $y_1$. We will use it in our relations for this section as well as in the following reduction:

Note that since the order of $|yx^z|=|xy^z|$ is $6$ we have, by properties of dihedral groups, that $(y^z)^{xy^z} \in C_G(x)$. We denote that element $y_2$. If the group $$C=\langle x,y,y_1,y_2\rangle$$ has index $1$ in $G$ then $G$ shrinks, since $C$ is clearly a subgroup of $C_G(x)$ and it having index 1 would imply that $z$ commutes with $x$.

Now akin to the previous section we will iterate through $G(6,i_2,i_3,i_4,i_5)$, where

$$w_1=xy^z,\ w_2=xyz,\ w_3=yy_1, w_4=xz^{yz},\ w_5=xx^{zyzxz}.$$

The above reduction allows us to discover most of the shrinking presentations. Some however require additional relations: in order to compute the index $|G:C|$ for the group $G(6,10,5,6,6)$ we need to additionally restrict the order of $|x^zx y^zy|$, while for the groups $G(6,12,\ldots)$ adding relations for $y^{zy} x^{zx}$ and $zy^{zxyzxzy}$ will enable one to proceed in a similar fashion.

We now present the remaining non-shrinking cases. 

The group $G(6,8,4,6)/\langle(y^zx^{(zxy)})^5\rangle$  has order $2^7 \cdot 3^3 \cdot 5 \cdot 11$. 
It is isomorphic to $2 \times M_{12}$ with
$$x=0 \times (1,11)(2,9)(3,12)(4,5)(6,10)(7,8),$$
$$y=1 \times (1,5)(2,7)(3,12)(4,11)(6,10)(8,9),$$
$$z=0 \times (1,6)(2,8)(3,7)(4,5)(9,11)(10,12).$$

Another non-shrinking case is $G(6,10,4,6,6)$. We compute that $|G(6,10,4,6,6):\langle z,y,y_1 \rangle|=660$. By constructing the action of that group on cosets we obtain that it is isomorphic to $(2.M_{22}):2$, meaning we are looking at the permutation representation of this group on $660$ points. 

Adding relations $(yz^{xz})^6$ and $(x^zy^{zyx})^6$ allows to determine the index of groups $\langle x,y_1,z \rangle$ and $\langle z^y, y^z, z^x \rangle$ with both of them having index $2464$ and yielding the aforementioned group, when constructing the action on cosets. The advantage of doing that is that these groups are fitting the case $(6,5)$ (two of the generators commute) so this reinforces the faithfulness of our presentation.

Group $G(6,10,6,5)/\langle(zy^{(zxyzxzy)})^{12},\ (xyy^z)^6\rangle$ has order $2^2\cdot 3^8 \cdot 5$. It is isomorphic to $O_3(G):D_{20}$.

Lastly, the group $G(6,12,4,4,4)$ is isomorphic to $2 \times (2^5:S_6)$, where the action of $S_6$ on the normal subgroup isomorphic to $2^5$ is again non-trivial, with $(5,6)$ and $(1,2,3,4,5)$ acting as matrices
$$
    \left(\begin{tabular}{c c c c c}
         0&0&0&0&1 \\
         1&1&0&0&1 \\
          0&0&0&1&1 \\
           1&0&1&0&0 \\
           1&0&0&0&0 \\
    \end{tabular}\right) \text{ and } \;
    \left(\begin{tabular}{c c c c c}
         1&0&0&1&0 \\
         1&0&1&1&0 \\
         0&0&0&0&1 \\
         0&1&0&1&0 \\
         1&0&0&0&0 \\
    \end{tabular}\right)
    $$
respectively.

\subsection{Additional remarks}

In this section we collect some interesting presentations for the case $r_1=r_2=6$.
Again $G(i_1\ldots i_7)$ denotes the group given by relations from the Diagram~\ref{InitDiag} and 
$$\{(xy^z)^{i_1},\ (xz^{yz})^{i_2},\ (xyz)^{i_3},\ (zz^{yzx})^{i_4},\ (xy^{zxyz})^{i_5},\ (zx^{zyx})^{i_6},\ (zy^{zyx}))^{i_7}\}.$$

\textbf{Examples} are as follows:\newline
$G(4,6,12,6,4,4,4) \simeq 2^6:(2^4:S_3 \times S_3)$,
$G(5,6,8,5,4) \simeq (SL(2,9):A_6):2^2$,\newline 
$G(5,6,10,6,6,5) \simeq 2^{10}:(2 \times PSL(2,11)$, 
$G(6,4,10,6,4,5,6) \simeq 2 \times (2.M_{22}:2).$\newline 
Moreover, $A_{12}\simeq G(6,6,12,6,4,5,5) / \langle (y^zy^{zyx})^4 \rangle$ with \newline
$x\simeq (1,2)(3,4)(5,6)(7,8)(9,10)(11,12),$ 
$y\simeq (1,2)(3,4)(5,7)(6,8)(9,12)(10,11),$ and \newline
$z\simeq (1,12)(2,5)(3,7)(4,6)(8,10)(9,11).$\newline
We note that $A_{12}$ is a maximal subgroup in the Monster group and so it had been studied extensively from the Majorana-theoretic point of view \cite{A12Maj}.

\section{Appendix: Presentations of groups}

In this section we provide the non-shrinking presentations that define maximal encountered $6$-transposition groups. 
All relations were stated over the course of the paper and this section exists for the reader's convenience.

\begin{center}
    
\begin{tabular}{ |l|c|c|c|c|}
\hline
N&Group & $(r_1,r_2) $&Relators& $|D|$\\ \hline
1&$2 \times S_4 \simeq 2^2:D_{12}$ & $(3,4)$ &  & $9$\\ \hline
2& $2 \times A_5$ & $(3,5)$& & $15$\\ \hline
3& $k^2:D_{12}$ & $(3,6)$& $(xyz)^{2k},\ k=4,5,6$ &$3\cdot 2k$\\ \hline
4& $k^2:D_8$ &$(4,4)$ & $(xy^z)^k,\ k=4,5,6$ & $4k$\\ \hline
5& $2^5:D_{10}$ & $(4,5)$ & $(xy^z)^4$, $(xz^{yz})^4$& $25$\\ \hline
6a& $PGL(2,9)$ & $(4,5)$ & $(xy^z)^5,\ (xyz)^8$ & $36$\\
6b& & $(5,5)$& $(xy^z)^4$ & \\ \hline
7a& $2 \times (2^4:S_5)$ &$(4,5)$& $(xy^z)^6,\  (xyz)^{12}$& $80$ \\ 
7b& &$(5,6)$ & $(xy^z)^4$, $(xyx^z)^{12}$ &\\
\hline
8& $2^6:D_{12}$ & $(4,6)$ & $(xy^z)^4$, $(xz^{yz})^4$ & $42$\\ \hline
9& $2^2 \times S_4$ & $(4,6)$ & $(xy^z)^4, (xz^{yz})^6, (xyz)^{12}$ & $15$\\ \hline 
10& $(A_5 \times A_5):2^2$ & $(4,6)$ & $(xy^z)^5,\  (xyz)^{12},\  (xz^{yz})^6$ &$160$ \\ \hline
11& $3^4:(D_8 \times S_3)$ & $(4,6)$ & $(xy^z)^6,\ (xyz)^{12},\  (xz^{yz})^4,\   (xy^zy)^6$ & $99$ \\ \hline
12& $(S_4 \times S_4):2^2$ & $(4,6)$ & $(xy^z)^6,\ (xyz)^{8},\  (xz^{yz})^4,$\  & $72$ \\ \hline
13& $2 \times (2^6:S_5)$ & $(4,6)$ & $(xy^z)^6,\  (xyz)^{10},\ (xz^{yz})^6,\  (z^yy^{zx})^4$ & $140$\\
\hline
14& $2^5:(2^4:D_{10})$ & $(5,5)$ & $(xy^z)^5,\ (xyz)^8$ &$80$ \\ \hline
15& $3^5:(2^4:D_{10})$ & $(5,5)$ & $(xy^z)^5,\ (xyz)^{12}$ & $180$\\ \hline
16a& $2^{10}:(2 \times PSL(2,11))$&$(5,5)$ & $(xy^z)^6,\ (xyz)^{10},\ (xz^{(xyz)^2})^4$ & $880$ \\ 
16b& & $(5,6)$ & $(xy^z)^5,\ (xyx^z)^{10},\ (x{x^z}^{(xyx^z)^2})^4$  & \\
\hline
17a& $O_2(G):A_5$& $(5,5)$ & $(xy^z)^6,\ (xyz)^{12},\ (xz^{xyz})^4$ & $240$ \\ 
17b& & $(5,6)$& $(xy^z)^5,\ (xyx^z)^{12},\ (z^yx^z)^4$ & \\\hline
18& $2 \times M_{12}$ & $(5,6)$ & \vtop{\hbox{\strut$(xy^z)^6,\ (xyz)^8,\ (yy^{zy(xz)^2})^4,$}
\hbox{\strut $(xz^{yz})^6,\  (y^zx^{zxy})^5$}}
& $396$ \\ \hline
19& $(2.M_{22}):2$ & $(5,6)$& \vtop{\hbox{\strut$(xy^z)^6,\ (xyz)^{10},\ (yy^{zy(xz)^2})^4,$\ }
\hbox{\strut $(xz^{yz})^6,\ (xx^{zyzxz})^6.\ H=\langle y, z, y^{zy(xz)^2}\rangle.$}} & $1485$ \\ \hline 

20& $O_3(G):D_{20}$ & $(5,6)$ &  \vtop{\hbox{\strut$(xy^z)^6,\ (xyz)^{10},\ (yy^{zy(xz)^2})^6,$\ } \hbox{\strut $(xz^{yz})^5,\ (zy^{zxyzxzy})^6,\ (xyy^z)^6$}} & $486$\\ \hline
21& $2 \times 2^5:S_6$ & $(5, 6)$ & $(xy^z)^6, (xyz)^{12}, (yy^{zy(xz)^2})^4, (xz^{yz})^4,\ (xx^{zyzxz})^4$ & $240$
\\ \hline
22& $2 \times 3.S_6$ & $(5,6)$ & $(xy^z)^6,\ (xyz)^{12},\ (yy^{zy(xz)^2})^6,\ (xz^{yz})^4$& $90$\\ \hline

\end{tabular}

\end{center}

We remark that some groups in the table are related, namely $1)$ and $2)$ are homomorphic images of $9)$ and $17)$ respectively.

\end{document}